\documentclass[12pt]{article}
\usepackage{latexsym,amsfonts,amsmath,amssymb}
\newtheorem{theorem}{Theorem}
\newtheorem{corollary}[theorem]{Corollary}
\newtheorem{sublemma}{Lemma}[theorem]
\newtheorem{lemma}[theorem]{Lemma}
\newtheorem{question}[theorem]{Question}
\newtheorem{observation}[theorem]{Observation}
\newtheorem{claim}[theorem]{Claim}
\newtheorem{subclaim}{Claim}[sublemma]
\newtheorem{conjecture}[theorem]{Conjecture}
\newtheorem{fact}[theorem]{Fact}
\newtheorem{definition}[theorem]{Definition}
\newtheorem{remark}[theorem]{Remark}
\newtheorem{example}[theorem]{Example}
\newtheorem{exercise}{Exercise}[section]
\def\Theorem #1.#2 #3\par{\def\claimname{Theorem}\setbox1=\hbox{#1}\ifdim\wd1=0pt
   \begin{theorem}{\rm #2} #3\end{theorem}\else
   \newtheorem{#1}[theorem]{#1}\begin{#1}\label{#1}{\rm #2} #3\end{#1}\fi}
\def\Corollary #1.#2 #3\par{\def\claimname{Corollary}\setbox1=\hbox{#1}\ifdim\wd1=0pt
   \begin{corollary}{\rm #2} #3\end{corollary}\else
   \newtheorem{#1}[theorem]{#1}\begin{#1}\label{#1}{\rm #2} #3\end{#1}\fi}
\def\Lemma #1.#2 #3\par{\def\claimname{Lemma}\setbox1=\hbox{#1}\ifdim\wd1=0pt
   \begin{lemma}{\rm #2} #3\end{lemma}\else
   \newtheorem{#1}[theorem]{#1}\begin{#1}\label{#1}{\rm #2} #3\end{#1}\fi}
\def\SubLemma #1.#2 #3\par{\def\claimname{Lemma}\setbox1=\hbox{#1}\ifdim\wd1=0pt
   \begin{sublemma}{\rm #2} #3\end{sublemma}\else
   \newtheorem{#1}{#1}[theorem]\begin{#1}\label{#1}{\rm #2} #3\end{#1}\fi}
\def\Question #1.#2 #3\par{\def\claimname{Question}\setbox1=\hbox{#1}\ifdim\wd1=0pt
   \begin{question}{\rm #2} #3\end{question}\else
   \newtheorem{#1}[theorem]{#1}\begin{#1}\label{#1}{\rm #2} #3\end{#1}\fi}
\def\Observation #1.#2 #3\par{\def\claimname{Observation}\setbox1=\hbox{#1}\ifdim\wd1=0pt
   \begin{observation}{\rm #2} #3\end{observation}\else
   \newtheorem{#1}[theorem]{#1}\begin{#1}\label{#1}{\rm #2} #3\end{#1}\fi}
\def\Claim #1.#2 #3\par{\def\claimname{Claim}\setbox1=\hbox{#1}\ifdim\wd1=0pt
   \begin{claim}{\rm #2} #3\end{claim}\else
   \newtheorem{#1}[theorem]{#1}\begin{#1}\label{#1}{\rm #2} #3\end{#1}\fi}
\def\SubClaim #1.#2 #3\par{\def\claimname{Claim}\setbox1=\hbox{#1}\ifdim\wd1=0pt
   \begin{subclaim}{\rm #2} #3\end{subclaim}\else
   \newtheorem{#1}{#1}[sublemma]\begin{#1}\label{#1}{\rm #2} #3\end{#1}\fi}
\def\Conjecture #1.#2 #3\par{\def\claimname{Conjecture}\setbox1=\hbox{#1}\ifdim\wd1=0pt
   \begin{conjecture}{\rm #2} #3\end{conjecture}\else
   \newtheorem{#1}[theorem]{#1}\begin{#1}\label{#1}{\rm #2} #3\end{#1}\fi}
\def\Fact #1.#2 #3\par{\def\claimname{Fact}\setbox1=\hbox{#1}\ifdim\wd1=0pt
   \begin{fact}{\rm #2} #3\end{fact}\else
   \newtheorem{#1}[theorem]{#1}\begin{#1}\label{#1}{\rm #2} #3\end{#1}\fi}
\def\Definition #1.#2 #3\par{\def\claimname{Definition}\setbox1=\hbox{#1}\ifdim\wd1=0pt
   \begin{definition}{\rm #2} {\rm #3}\end{definition}\else
   \newtheorem{#1}[theorem]{#1}\begin{#1}\label{#1}{\rm #2} {\rm #3}\end{#1}\fi}
\def\Remark #1.#2 #3\par{\def\claimname{Remark}\setbox1=\hbox{#1}\ifdim\wd1=0pt
   \begin{remark}{\rm #2} {\rm #3}\end{remark}\else
   \newtheorem{#1}[theorem]{#1}\begin{#1}\label{#1}{\rm #2} {\rm #3}\end{#1}\fi}
\def\Example #1.#2 #3\par{\def\claimname{Example}\setbox1=\hbox{#1}\ifdim\wd1=0pt
   \begin{example}{\rm #2} #3\end{example}\else
   \newtheorem{#1}[theorem]{#1}\begin{#1}\label{#1}{\rm #2} #3\end{#1}\fi}
\def\Exercise #1.#2 #3\par{\def\claimname{Exercise}\setbox1=\hbox{#1}\ifdim\wd1=0pt
   {\footnotesize\begin{exercise}{\rm #2} {\rm #3}\end{exercise}}\else
   \newtheorem{#1}[section]{#1}{\footnotesize\begin{#1}\label{#1}{\rm #2} {\rm #3}\end{#1}}\fi}
\def\QuietTheorem #1.#2 #3\par{\setbox1=\hbox{#1}\ifdim\wd1=0pt\proclaim{Theorem {\rm #2}}{#3}\else\proclaim{#1 {\rm #2}}{#3}\fi}
\newcommand{\proclaim}[2]{\smallskip\noindent{\bf #1} {\sl#2}\par\smallskip}
\def\Proclaim #1.#2 #3\par{\proclaim{#1 {\rm #2}}{#3}}
\newenvironment{proof}{\noindent}{\kern2pt\QEDbox\par\medskip}
\newif\ifinproof\inprooffalse
\def\Proof#1: {\setbox1=\hbox{#1}\ifdim\wd1=0pt\begin{proof}{\bf Proof: }\else\medskip\begin{proof}{\bf #1: }\fi\inprooftrue}
\newcommand{\QED}{\end{proof}}

\def\BF#1.{{\bf #1.}}
\def\Abstract #1\par{\begin{quotation}{\singlespaced\footnotesize{\noindent{\bf Abstract.~}#1}}\end{quotation}}
\def\Title #1\par{\title{#1}\maketitle}
\def\Author #1\par{\author{#1}}
\def\Acknowledgement#1\par{\thanks{#1}}
\def\Chapter #1\par{\chapter{#1}}
\def\Section #1\par{\section{#1}}
\def\QuietSection #1\par{\section*{#1}}
\def\SubSection #1\par{\subsection{#1}}
\def\SubSubSection #1\par{\subsubsection{#1}}
\def\MidTitle #1\par{\bigskip\goodbreak\centerline{\small\bf #1}\bigskip\noindent}
\def\Margin #1\par{\marginpar{\tiny #1}}

\newcommand{\singlespaced}{\baselineskip=15pt}
\newcommand{\citindex}{\index[citdep]}
\newcommand{\Ref}[2][Theorem]{\edef\entry{\ref{#2}..@#1 \ref{#2}!\ifinproof\arabic{chapter}.\arabic{section}.\arabic{theorem}..@\claimname\ \arabic{chapter}.\arabic{section}.\arabic{theorem}\else \arabic{chapter}.\arabic{section}..@Section \arabic{chapter}.\arabic{section}\fi|dotfill\ p.}\expandafter\citindex\expandafter{\entry}\ref{#2}}
\def\bottomnote #1\par{{\renewcommand{\thefootnote}{}\footnotetext{#1}}}
\newcommand{\N}{{\mathbb N}}

\newfont{\msam}{msam10 at 12pt}

\newcommand{\set}[1]{\{\,{#1}\,\}}

\newcommand{\cross}{\times}

\newcommand{\union}{\cup}

\newcommand{\intersect}{\cap}

\newcommand{\smalllt}{\mathrel{\mathchoice{\raise2pt\hbox{$\scriptstyle<$}}{\raise1pt\hbox{$\scriptstyle<$}}{\raise0pt\hbox{$\scriptscriptstyle<$}}{\scriptscriptstyle<}}}
\newcommand{\smallleq}{\mathrel{\mathchoice{\raise2pt\hbox{$\scriptstyle\leq$}}{\raise1pt\hbox{$\scriptstyle\leq$}}{\raise1pt\hbox{$\scriptscriptstyle\leq$}}{\scriptscriptstyle\leq}}}

\newcommand{\boolval}[1]{\mathopen{\lbrack\!\lbrack}\,#1\,\mathclose{\rbrack\!\rbrack}}
\def\[#1]{\boolval{#1}}

\newcommand{\UnderTilde}[1]{{\setbox1=\hbox{$#1$}\baselineskip=0pt\vtop{\hbox{$#1$}\hbox to\wd1{\hfil$\sim$\hfil}}}{}}
\newcommand{\Undertilde}[1]{{\setbox1=\hbox{$#1$}\baselineskip=0pt\vtop{\hbox{$#1$}\hbox to\wd1{\hfil$\scriptstyle\sim$\hfil}}}{}}
\newcommand{\undertilde}[1]{{\setbox1=\hbox{$#1$}\baselineskip=0pt\vtop{\hbox{$#1$}\hbox to\wd1{\hfil$\scriptscriptstyle\sim$\hfil}}}{}}
\newcommand{\UnderdTilde}[1]{{\setbox1=\hbox{$#1$}\baselineskip=0pt\vtop{\hbox{$#1$}\hbox to\wd1{\hfil$\approx$\hfil}}}{}}
\newcommand{\Underdtilde}[1]{{\setbox1=\hbox{$#1$}\baselineskip=0pt\vtop{\hbox{$#1$}\hbox to\wd1{\hfil\scriptsize$\approx$\hfil}}}{}}

\renewcommand{\th}{{\hbox{\scriptsize th}}}

\def\<#1>{\langle#1\rangle}

\newcommand{\QEDbox}{\fbox{}}

\newcommand{\cell}[1]{\boxit{\hbox to 17pt{\strut\hfil$#1$\hfil}}}
\newcommand{\head}[2]{\lower2pt\vbox{\hbox{\strut\footnotesize\it\hskip3pt#2}\boxit{\cell#1}}}
\newcommand{\boxit}[1]{\setbox4=\hbox{\kern2pt#1\kern2pt}\hbox{\vrule\vbox{\hrule\kern2pt\box4\kern2pt\hrule}\vrule}}
\newcommand{\Col}[3]{\hbox{\vbox{\baselineskip=0pt\parskip=0pt\cell#1\cell#2\cell#3}}}
\newcommand{\tapenames}{\raise 5pt\vbox to .7in{\hbox to .8in{\it\hfill input: \strut}\vfill\hbox to
.8in{\it\hfill scratch: \strut}\vfill\hbox to .8in{\it\hfill output: \strut}}}
\newcommand{\Head}[4]{\lower2pt\vbox{\hbox to25pt{\strut\footnotesize\it\hfill#4\hfill}\boxit{\Col#1#2#3}}}
\newcommand{\Dots}{\raise 5pt\vbox to .7in{\hbox{\ $\cdots$\strut}\vfill\hbox{\ $\cdots$\strut}\vfill\hbox{\
$\cdots$\strut}}}
\renewcommand{\dots}{\raise5pt\hbox{\ $\cdots$}}
\newcommand{\factordiagramup}[6]{$$\begin{array}{ccc}
#1&\raise3pt\vbox{\hbox to60pt{\hfill$\scriptstyle
#2$\hfill}\vskip-6pt\hbox{$\vector(4,0){60}$}}&#3\\ \vbox
to30pt{}&\raise22pt\vtop{\hbox{$\vector(4,-3){60}$}\vskip-22pt\hbox
to60pt{\hfill$\scriptstyle #4\qquad$\hfill}}
     &\ \ \lower22pt\hbox{$\vector(0,3){45}$}\ {\scriptstyle #5}\\
\vbox to15pt{}&&#6\\
\end{array}$$}
\newcommand{\factordiagram}[6]{$$\begin{array}{ccc}
#1&&\\ \ \ \raise22pt\hbox{$\vector(0,-3){45}$}\ {\scriptstyle #2}
&\raise22pt\hbox{$\vector(2,-1){90}$}\raise5pt\llap{$\scriptstyle#3$\qquad\quad}&\vbox
to25pt{}\\ #4&\raise3pt\vbox{\hbox to90pt{\hfill$\scriptstyle
#5$\hfill}\vskip-6pt\hbox{$\vector(4,0){90}$}}&#6\\
\end{array}$$}
\newcommand{\df}{\it} 
\begin{document}
\author{Joel David Hamkins\\
\normalsize\sc The City University of New York\\
{\footnotesize http://jdh.hamkins.org}\\
\\
Alexei Miasnikov\\
\normalsize\sc The City University of New York\\
{\footnotesize  http://www.cs.gc.cuny.edu/$\sim$amyasnikov}}

\bottomnote MSC: 03D10; 68Q17. Keywords: Turing machines, halting problem, decidability. The first author is affiliated
with the College of Staten Island of CUNY and The CUNY Graduate Center, and his research has been supported by grants
from the Research Foundation of CUNY. The second author is affiliated with The City College of New York and The CUNY
Graduate Center.

\Title The halting problem is decidable on a set of asymptotic probability one

\Abstract The halting problem for Turing machines is decidable on a set of asymptotic probability one. The proof is
sensitive to the particular computational model.

The halting problem for Turing machines is perhaps the canonical undecidable set. Nevertheless, we prove that there is
an algorithm deciding almost all instances of it. The halting problem is therefore among the growing collection of
those exhibiting the ``black hole'' phenomenon of complexity theory, by which the difficulty of an unfeasible or
undecidable problem is confined to a very small region, a black hole, outside of which the problem is easy.

We use the most natural method for measuring the size of a set of Turing machine programs, namely, that of {\df
asymptotic density}. The asymptotic density or probability of a set $B$ of Turing machine programs is the limit of the
proportion of all $n$-state programs in $B$ as $n$ increases. That is, if $P_n$ is the set of all $n$-state programs,
then the asymptotic probability of $B$ is
$$\mu(B)=\lim_{n\to\infty}{|B\intersect P_n|\over|P_n|},$$ provided that this limit exists. If $B$ has asymptotic probability one, for example, then for
sufficiently large $n$, more than 99\% of all $n$-state programs are in $B$, and so on as close to 100\% as desired.

\Theorem Main Theorem. There is a set $B$ of Turing machine programs such that
\begin{enumerate}
 \item $B$ has asymptotic probability one.
 \item $B$ is polynomial time decidable.
 \item The halting problem $H\intersect B$ is polynomial time decidable.
\end{enumerate}

The proof is sensitive to the particular (but common) computational model. We use the Turing machine model with a finite program directing the
operation of a head reading and writing $0$s and $1$s while moving on a one-way infinite
tape.$$\cell1\cell0\cell1\lower2pt\vbox{\hbox{\strut\hskip2pt\small Head}\boxit{\cell1}}\cell1\cell0\cell0\raise5pt\hbox{\ $\cdots\longrightarrow$}$$
The Turing machine has $n$ states $Q=\set{q_1,\ldots,q_n}$, with $q_1$ designated as the {\it start} state, plus a separate designated {\it halt}
state, which is not counted as one of the $n$ states. A Turing machine {\it program} is a function
$$p:Q\cross\set{0,1}\to (Q\union\set{{\it halt}})\cross\set{0,1}\cross\set{L,R}.$$ The transition $p(q,i)=\<r,j,R>$,
for example, directs that when the head is in state $q$ reading symbol $i$, it should change to state $r$, write symbol $j$, and move one cell to the
right. The computation of a program proceeds by iteratively performing the instructions of such transition rules, halting when (and if) the {\it halt}
state is reached. If the machine attempts to move left from the left-most cell, then the head falls off the tape and all computation ceases. Since the
domain of the program has size $2n$ and the target space has size $4(n+1)$, we can easily count the number of programs:

\SubLemma. The number of $n$-state Turing machine programs is $(4(n+1))^{2n}$.

As a warm-up exercise, let us calculate the asymptotic probability of the set of programs having no transition reaching the {\it halt} state. Such a
criterion is clearly linear time decidable (for any reasonable representation of programs by finite binary sequences), and no computation by such a
program can ever reach the {\it halt} state.

\SubLemma. The collection of programs having no transition reaching the {\it halt} state has asymptotic probability
$1/e^2$, which is about $13.5\%$.\label{NoHaltState}

\Proof: If $p$ has no transition reaching the {\it halt} state, then $p:Q\cross\set{0,1}\to
Q\cross\set{0,1}\cross\set{L,R}$. Since this target set has size $4n$, the total number of such functions is
$(4n)^{2n}$. The asymptotic proportion of all $n$-state programs with this property is therefore
$$\lim_{n\to\infty}{ (4n)^{2n}\over (4(n+1))^{2n}}=\lim_{n\to\infty}\left({n\over
n+1}\right)^{2n} =\lim_{n\to\infty}{\left[\bigl(1-{1\over n+1}\bigr)^n\right]^2}=1/e^2.$$ Therefore, the asymptotic
probability that a Turing machine program does not engage the {\it halt} state is $1/e^2$.\QED

\Definition. The {\df halting problem} is the set $H$ of programs $p$ that halt when computing on input $0$, on a tape
initially filled with $0$s.

For the purposes of defining the halting problem $H$, one should  specify whether it officially counts as halting or
not, if the head should happen to fall off the left edge of the tape. For definiteness, let us regard such computations
as having not officially halted, as the {\it halt} state was not reached. Thus, we regard $H$ as the set of programs
that eventually reach the {\it halt} state from an initially empty tape. To be even more specific, if the head happens
to fall off the tape while executing the transition $p(q,i)=\<r,j,L>$, then we do not regard the state $r$ as having
been achieved, since this step was not completed.

\Proof Proof of Main Theorem: We now prove the Main Theorem. Let $B$ be the set of programs that on input $0$ either
halt before repeating a state or fall off the tape before repeating a state. Clearly, $B$ is polynomial time decidable,
since we need only run a program $p$ for at most $n$ steps, where $n$ is the number of states in $p$, to determine
whether or not it is in $B$. It is equally clear that the halting problem is polynomial time decidable for programs $p$
in $B$, since again we need only simulate $p$ for $n$ steps to know whether it halted or fell off. What remains is to
prove that this behavior occurs with asymptotic probability one.

\addtocounter{theorem}{-1}\addtocounter{sublemma}{2}

\SubLemma. For any fixed input and fixed $k\geq 0$, the set of programs not repeating states within the first $k$ steps
 has asymptotic probability one.\label{NoRepeats}

\Proof: Just to be clear, we count a computation that halts or falls off the tape as satisfying the property, if it
does so before repeating a state. We calculate for large $n$ the proportion of all $n$-state programs having this
property, by induction on $k$. When $k=0$, then all programs have the property. Suppose that the set $B_k$ of programs
having the desired property for $k$ has asymptotic probability one, and consider $B_{k+1}$. Fix any $\epsilon$, and
choose $n$ large enough so that $B_k$ has proportion more than $1-\epsilon/2$ of all $n$-state programs. Among all
$n$-state programs $p$ in $B_k$, consider the probability that $p$ is in $B_{k+1}$. If $p$ leads to a computation where
the head has already fallen off the tape, then of course $p\in B_{k+1}$. Otherwise, the first $k$ steps of computation
by $p$ have led to the successive states $q_{i_0}, q_{i_1}, \ldots,q_{i_k}$, which have not yet repeated. The
$(k+1)^\th$ step of computation involves a transition rule $p(q_{i_k},j_k)=(q_{i_{k+1}},j_{k+1},m_k)$, giving
respectively the new state, the new bit to write on the tape and the direction to move the head. In order for $p$ to be
in $B_{k+1}$, it suffices that $q_{i_{k+1}}$ must not be one of the previously used states $\set{q_{i_0}, q_{i_1},
\ldots,q_{i_k}}$. Since there are $(n+1)-(k+1)=n-k$ many other equally likely states to choose from, the proportion of
all $n$-state programs agreeing with $p$ on the first $k$ steps and satisfying the additional requirement is $n-k\over
n$. The proportion of all $n$-state programs in $B_{k+1}$, consequently, is at least $(1-\epsilon/2)({n-k\over n})$.
Since $n-k\over n$ goes to $1$ as $n$ becomes large, we may choose $n$ large enough so that this proportion is at least
$1-\epsilon$. Thus, $B_{k+1}$ has asymptotic probability one, as desired.\QED

Proceeding with the main argument, let $B_k$ be the set of programs that do not repeat a state within their first $k$
steps of computation. The key idea is that for the first $k$ steps of computation, the programs in $B_k$ behave
statistically like a random walk with uniform probability of going left or right. The reason is that if a program lands
in a totally new state $q$, reading some symbol $i$, then among the programs landing in that situation and agreeing
with the computation so far, exactly half of them will opt to move left and half will move right, precisely because
nothing about state $q$ has yet been determined. Because of this, we may make use of Polya's classical result on random
walks, which we mention without proof.

\SubLemma.({Polya \cite{Polya1921:IrrfahrtImStrassennetz}, see also e.g. \cite{Feller1968:IntroProbabilityTheory}}) In
the random walk with equal likelihood of moving left or right on a one-way infinite tape, beginning on the left-most
cell, the probability of eventually falling off the left edge is $1$.\label{RandomWalk}

This is the famous {\it recurrence} phenomenon, because it asserts that such a random walk has probability one of
eventually returning to its starting point. It follows that with probability one the random walk reaches any given
fixed position of the tape. Interestingly, the recurrence property holds for random walks in dimensions one and two,
but not in dimensions three or higher.

Putting everything together, let us show that $B$ has asymptotic probability one. Fix any $\epsilon>0$, and by Lemma
\ref{RandomWalk} find some large $k$ such that with probability exceeding $\sqrt{1-\epsilon}$, the $k$-step random walk
falls off the left edge of the tape. By Lemma \ref{NoRepeats}, let $n$ be large enough so that $B_k$ contains more than
the proportion $\sqrt{1-\epsilon}$ of all $n$-state programs. Combining these facts with the observation that programs
in $B_k$ operate statistically like random walks for their first $k$ steps of computation (or until they halt, if this
is sooner), as far as the head position is concerned, we conclude that proportion at least
$(\sqrt{1-\epsilon})^2=1-\epsilon$ of all $n$-state programs exhibit the desired property. So the set $B$ of all such
programs has asymptotic probability one, and the theorem is proved.\QED

\addtocounter{theorem}{1}

Let us now clarify matters by untangling the two possibilities for programs in $B$, namely, (1) the programs that halt before repeating a state and (2)
the programs that fall off the tape before repeating a state. The fact is that behavior (1) is very rare and behavior (2) occurs with asymptotic
probability one.

\Theorem. The asymptotic probability one behavior of a Turing machine, on any fixed input, is that the head falls off the tape before halting or
repeating a state.\label{Prob1Behavior}

\Proof: First, we generalize Lemma \ref{NoRepeats} to exclude the possibility of halting.

\SubLemma. For any fixed input and fixed $k\geq 0$, the set of programs not repeating states and not halting within the first $k$ steps of computation
on that input has asymptotic probability one.\label{NoRepeatsOrHalts}

\Proof: Let $C_k$ be the desired set of programs, which includes the programs that fall off the tape within the first $k$ steps of computation on that
fixed input, provided that they do so before repeating a state. As in Lemma \ref{NoRepeats}, we show inductively that $C_k$ has asymptotic density one.
When $k=0$, this is trivial. Let us now calculate the probability that a program $p$ is in $C_{k+1}$, given that it is in $C_k$. If the head fell off
within $k$ steps, then $p$ will also be in $C_{k+1}$. Otherwise, as in Lemma \ref{NoRepeats}, the first $k$ steps of computation exhibit states
$q_{i_0}, q_{i_1}, \ldots,q_{i_k}$, which have not yet repeated. The $(k+1)^\th$ step of computation involves a transition
$p(q_{i_k},j_k)=(q_{i_{k+1}},j_{k+1},m_k)$, which will place $p$ into $C_{k+1}$ if $q_{i_{k+1}}$ is a new state and not the {\it halt} state. Since
there are $n-(k+1)$ remaining states to choose from, the probability that $p$ will be in $C_{k+1}$ is at least $n-(k+1)\over n$. Since this probability
goes to $1$ as $n$ goes to infinity, we conclude that $C_{k+1}$ has asymptotic probability one.\QED

Now fix any $\epsilon>0$. Select $k$ large enough so that the random walk in $k$ steps has probability exceeding
$\sqrt{1-\epsilon}$ of falling off the left edge. By Lemma \ref{NoRepeatsOrHalts}, take $n$ sufficiently large so that
the proportion of all $n$-state programs that do not halt in $k$ steps and do not repeat a state in $k$-steps is at
least $\sqrt{1-\epsilon}$. Thus, as in the Main Theorem, these computations behave statistically like random walks, as
far as the head position is concerned, and so the proportion of all $n$-state machines that fall off the tape in $k$
steps before repeating a state or halting is at least $\sqrt{1-\epsilon}\sqrt{1-\epsilon}=1-\epsilon$, as desired.\QED

\Corollary. The halting problem $H$ has asymptotic probability zero. And the complement of $H$ contains a decidable set of asymptotic probability
one.\label{HisDensity0}

\Proof: If the head falls off the tape, then the computation cannot reach the {\it halt} state, and so the program is not in $H$. So $H$ has density
zero. The set of programs that fall off the tape before repeating a state is contained in the complement of $H$, is clearly polynomial time decidable
and, by Theorem \ref{Prob1Behavior}, has asymptotic probability one.\QED

The previous Corollary depends on the formalism that computations for which the head falls off the tape are not counted as halting. If one wishes
instead to count them as halting, then the conclusion would be that the corresponding version of $H$ would have asymptotic probability one and contain
a decidable set of asymptotic probability one.

Because the computational behavior identified in Theorem \ref{Prob1Behavior}, with the head falling off the tape before a state is repeated, is both
typical and trivial, many other well-known undecidability problems for Turing machines can also be decided with asymptotic probability one. We give two
examples.

\Definition. Let FIN  be the set of programs computing functions on $\N$ with finite domain and COF be the set of programs with cofinite domain. These
sets are well known to be undecidable (see e.g. \cite{Soare1987:ReSetsAndDegrees}).

\Corollary. There is a set of programs $B$ such that:
\begin{enumerate}
 \item $B$ has asymptotic probability one.
 \item $B$ is polynomial time decidable.
 \item $FIN\intersect B$ is polynomial time decidable.
 \item $COF\intersect B$ is polynomial time decidable.
\end{enumerate}\label{FinInfCof}

\Proof: For the purposes of computing functions on $\N$, we assume that input on a Turing machine tape is given by a string of $1$s, that is, in unary
form. Let $B$ be the set of programs that fall off the tape before halting or repeating a state, on a tape initially filled entirely with $1$s. This
set is clearly polynomial time decidable, and by Theorem \ref{Prob1Behavior}, it has asymptotic probability one. But any program that falls off the
tape will have had a chance to inspect only finitely many of the $1$s on the tape before doing so, and so the program will have this same behavior
provided that there is a sufficiently long string of $1$s on the tape as input. So every program in $B$ is in FIN and none are in COF. On $B$,
therefore, these questions are decidable.\QED

The proof of Corollary \ref{FinInfCof} shows that almost every program computes a finite function. In other words, FIN
has asymptotic probability one and COF has asymptotic probability zero. If one takes the domains of the computable
functions as the natural enumeration of the computably enumerable (c.e.) sets, then this means that almost every c.e.
set is finite.

Let us turn now to the question of whether the conclusions of the main theorem hold for other models of computability.

\Corollary. The conclusion of the main theorem also holds for the following models of computability:
\begin{enumerate}
 \item Single tape Turing machines with an arbitrary finite alphabet, operating on a one-way infinite tape.
 \item Multi-tape Turing machines with an arbitrary finite alphabet, operating on one-way infinite tapes.
 \item Turing machines with a head moving on a half-plane or quarter-plane grid of cells with an arbitrary finite
        alphabet.
\end{enumerate}

\Proof: For the multi-tape model of (2), we assume that there is a single head moving back and forth, reading and
writing on all columns at once. The corollary is proved merely by observing that the calculations of Lemma
\ref{NoRepeats} do not fundamentally rely on the size of the alphabet, and so in the case of a general alphabet, it is
still true that for any $k$ the set of programs that adopt new states for their first $k$ moves has asymptotic
probability one. Because these programs therefore act like a random walk for the first $k$ steps, the probability that
they fall off the left end of the tape can be made as close to $1$ as possible. So (1) holds. The multi-tape model of
(2) is functionally equivalent to having a larger alphabet, if one regards a entire column of cell values as a single
element of a larger alphabet. So (2) holds. For (3), we observe that first, the analogue of Lemma \ref{NoRepeats}
remains true, and second, the desired conclusion now follows from the two-dimensional generalization of Lemma
\ref{RandomWalk}, by which random walks on a the half-plane (or any smaller portion of the plane), eventually fall off
the edge with probability one.\QED

The result also applies to 2-dimensional Turing machines operating on a full doubly-infinite plane, provided that it
has at least one broken cell, which is broken in the sense that it causes the computation to cease if the head should
happen to occupy it. The point is that because of the 2-dimensional analogue of Polya's recurrence theorem, on any
fixed input such a Turing machine would with asymptotic probability one land on the forbidden cell before repeating a
state.

The reader will have already observed, of course, that our argument does {\it not} work with Turing machines using
doubly-infinite tapes, another common model, for which there is no possibility that the head falls off the tape. And
neither does it work with the one-way infinite tape models that allow computation somehow to continue after attempting
to move left from the left-most cell. We admit that this situation is unsatisfactory, because one doesn't like results
in computability theory to be sensitive to the choice of computational model.

\Question. Does the conclusion of the Main Theorem hold for all models of Turing machine computation?

The focus, of course, is on the models with two-way infinite tapes. One can weaken the desired conclusion by asking
only that the halting problem be decided on a set of large probability, rather than probability one. If one weakens
this too much, by asking only to decide the problem on a set of nonzero probability, then it becomes a trivial
consequence of Lemma \ref{NoHaltState}:

\Theorem. For any model of Turing machine, including those with two-way infinite tapes, there is a set $B$ of Turing
machine programs such that
\begin{enumerate}
 \item $B$ has nonzero asymptotic probability.
 \item $B$ is polynomial time decidable.
 \item The halting problem $H\intersect B$ is polynomial time decidable.
\end{enumerate}\label{MostPrograms}

\Proof: The set of programs arising in Lemma \ref{NoHaltState}, which have no transition leading to the {\it halt}
state, has asymptotic probability $1/e^2$; but clearly no such program is in $H$.\QED

\Question. In Theorem \ref{MostPrograms}, how large can the probability of $B$ be? Can one always decide the halting
problem in an asymptotic majority of cases?

We close the paper with a few elementary cautionary observations. First, the relation of Turing equivalence does not
respect asymptotic density.

\Theorem. The relation of Turing equivalence does not respect the property of having asymptotic probability one in the
natural numbers. Indeed, for any set $A$ of natural numbers there is a set $B$ that is Turing equivalent to $A$ and has
any prescribed asymptotic density (or nonconvergent density).

\Proof: If a set $A$ has asymptotic density one in the natural numbers, then the complement of $A$, which is Turing
equivalent to $A$, has asymptotic density zero. But also, any set $A$ is Turing equivalent to a set with asymptotic
density zero, by simply multiplying its second member by $2$, it's third member by $3$, and so on, so as to stretch it
out to density zero. The complement of this set, which is also Turing equivalent, has asymptotic density one.
Intermediate densities can be achieved by adding regular blocks of numbers in a computable pattern, so as to achieve a
given intermediate density, while the true information is coded on a thin set of density zero. By alternating blocks of
numbers with large empty stretches, one can arrange that the asymptotic density of the set does not converge, and even
that the upper density is 1 while the lower density is 0. Meanwhile, the true information of the set is coded on a thin
set, of density zero, which does not upset those calculations.\QED

Second, the notion of what happens ``almost everywhere'' can be highly sensitive to what are otherwise unimportant
differences in formalism. For example, if one takes as the basic model of computability a suitable generalization of
C++ programs, then most would agree that for the usual purposes it is an irrelevant formalism whether one excludes
programs at the outset that have syntax errors preventing them from compiling, or instead takes them to compute the
empty function. But if they were officially counted, then because clearly there are far more programs with errors than
without, it would mean that almost every program would be trivial in this way. For such a model, all interesting
phenomena would occur on a set of asymptotic density zero.

\bibliographystyle{alpha}
\bibliography{MathBiblio,HamkinsBiblio}

\end{document}